\documentclass[12pt,leqno]{article}

\usepackage{latexsym,amsmath,amssymb,amsthm,amsfonts}

\renewcommand{\epsilon}{\varepsilon}

\newtheorem{theorem}{Theorem}

\begin{document}

\author{Andriy V.\ Bondarenko, Maryna S. Viazovska}
\title{Bernstein type inequality in monotone rational approximation}
\date{}
\maketitle

\begin{abstract}
The following analog of Bernstein inequality for monotone rational
functions is established: if $R$ is an increasing on $[-1,1]$
rational function of degree $n$, then
$$
R'(x)<\frac{9^n}{1-x^2}\|R\|,\quad x\in (-1,1).
$$
The exponential dependence of constant factor on $n$ is shown,
with sharp estimates for odd rational functions.
\end{abstract}

\noindent
Andriy V.\ Bondarenko, Maryna S.\ Viazovska\\ Faculty of Mech. and
Mathematics,\\ Kyiv National Taras Shevchenko university,\\ Kyiv,
01033, Ukraine.\\ tel: (+38)0442590591\\ e-mail: {\tt
bonda@univ.kiev.ua,  v-marina@ukr.net}

\newpage

\section {Introduction}
Let $P_n$ be the space of all polynomials of degree at most $n$.
Denote by $Q_n$ the set of all continuous rational functions on
$[\,-1,1\,]$, $r=\frac pq$, where $p$, $q\in P_n$. Now we state
the well-known Bernstein inequality: If $p\in P_n$, then
\begin{equation}
\label{i8}
|p'(x)|\le\frac n{\sqrt{1-x^2}}\,\|p\|,\quad x\in (-1,1),
\end{equation}
where $\|\cdot\|:=\|\cdot\|_{C[-1,1]}$ is the uniform norm on
$[-1,1]$. Unfortunately, the direct analog of this inequality in
rational approximation impossible to establish. Indeed, if
$R(x)=\frac{\delta\,x}{x^2+\delta^2}$, $\delta>0$, then $|R(x)|<
1$, $x\in \mathbb{R}$ and $R'(0)=\delta^{-1}$ can be arbitrary
large. What is true, is Pekarskii inequality [1], where the norm
of $R'$ and of $R$ are taken in different spaces. Our main result
is
\begin{theorem}
If $R\in Q_{2n}$ is an odd and increasing function on $[-1,1]$,
then
\begin{equation}
\label{i1}
R'(0)\le\frac 12\cdot 9^nR(1).
\end{equation}
\end{theorem}
Theorem 1 easily implies the following analog of
estimate~\eqref{i8} for all increasing on $[-1,1]$ rational
functions.\\ {\bf Corollary 1.} If $R\in Q_n$ is an increasing
function on $[-1,1]$, then
$$
R'(x)<\frac{9^n}{1-x^2}\|R\|,\quad x\in (-1,1).
$$
A lower estimate for the constant in the right hand side
of~\eqref{i1} is provided by
\begin{theorem}
For each $n\in\mathbb{N}$
$$
\sup_{R}\frac{R'(0)}{\|R\|}\ge 9^{n-1},
$$
where the supremum is taken over the set of all odd increasing on
$[-1,1]$ rational function $R\in Q_{2n-1}$.
\end{theorem}

In Section 2 we prove some auxiliary results, in Section 3 we
prove Theorem~1 and in Section 4 we prove Theorem 2.

\section{Auxiliary lemmas}
Let $u_i<v_i$, $i=\overline{1,n}$, be arbitrary numbers. Put
\begin{align*}
\Pi:&=\left\{\vec{y}=(y_1,\ldots,y_n)\in\mathbb{R}^n\,|\,u_i
\le y_i\le v_i, i=\overline{1,n}\right\},\\
\Pi^+_k:&=\left\{\vec{y}=(y_1,\ldots,y_n)\in\Pi\,|\,y_k=v_k\right\},
\quad k=\overline{1,n},
\end{align*}
and
$$
\Pi^-_k:=\left\{\vec{y}=(y_1,\ldots,y_n)\in\Pi\,|\,y_k=u_k\right\},
\quad k=\overline{1,n}.
$$
To prove the following Lemma 1 we use the well known Brouwer
fixed-point theorem \cite{N}\\
{\bf Theorem B.} Let A be a closed bounded convex subset of
$\mathbb{R}^n$ and $F:A\rightarrow A$ be a continuous mapping on
$A$. Then $F(\vec{z})=\vec{z}$, for some $\vec{z}\in A$.\\
{\bf Lemma 1.} Let $f_k:\Pi\rightarrow\mathbb{R}$,
$k=\overline{1,n}$, be continuous functions satisfying the
following inequalities
$$
f_k(\vec{y})<0,\quad \vec{y}\in\Pi_k^-,\quad k=\overline{1,n},
$$
and
$$
f_k(\vec{y})>0,\quad \vec{y}\in\Pi_k^+,\quad k=\overline{1,n}.
$$
Then, there exists $\vec{z}\in\Pi$ such that $f_k(\vec{z})=0$,
$k=\overline{1,n}$.
\begin{proof}
Without any loss of generality we may assume that $\Pi=[-1;1]^n$.
Put
$\varphi(\vec{x})=(\varphi_1(\vec{x}),...,\varphi_n(\vec{x}))$,
where
$$
\varphi_k(\vec{x})=\frac{f_k(\vec{x})}{|f_i(\vec{x})|+(1-x_k^2)}.
$$
This definition readily implies
\begin{gather}
\label{L1}
|\varphi_i(\vec{x})|\le 1,\qquad \vec{x}\in\Pi,\\
\label{L2}
\varphi_i(\vec{x})=1,\quad\vec{x}\in\Pi_k^+,\quad k=\overline{1,n},\\
\label{L3}
\varphi_i(\vec{x})=-1,\quad\vec{x}\in\Pi_k^-,\quad k=\overline{1,n}.
\end{gather}
Since each $\varphi_k$, $k=\overline{1,n}$, is a continuous
function on $\Pi$, then by~\eqref{L2} and~\eqref{L3}, there exists
a number $\mu >0$ small enough, such that
\begin{equation}
\label{L4}
\varphi(\vec{x})>0,\quad \vec{x}=(x_1,\ldots, x_n)\in\Pi,\quad 1-\mu\le x_k\le
1,
\end{equation}
and
\begin{equation}
\label{L5}
\varphi(\vec{x})<0,\quad \vec{x}=(x_1,\ldots, x_n)\in\Pi,\quad -1\le x_k\le
-1+\mu.
\end{equation}
Now we prove that the set $A:=\Pi$ and the mapping
$F(\vec{x}):=\vec{x}-\mu\varphi(\vec{x})$ satisfy the conditions
of the Theorem B. Since $\varphi$ is a continuous mapping, then
$F$ is a continuous mapping as well. Finally, we prove that
$F(\vec{x})\in\Pi$, for all $\vec{x}\in\Pi$, that is
\begin{equation}
\label{L6}
-1\le x_k-\mu\varphi_k(\vec{x})\le 1,\quad x_k\in [\,-1,1\,],\quad
k=\overline{1,n}.
\end{equation}
If $x_k\in [\,-1+\mu, 1-\mu\,]$, then the inequality~\eqref{L6}
readily follows from~\eqref{L1}. Taking into account~\eqref{L4}
and~\eqref{L5} we get~\eqref{L6} for $x_k\in [\,1-\mu,1\,]$ and
for $x_k\in [\,-1,-1+\mu\,]$ respectively, $k=\overline{1,n}$,
so~\eqref{L6} holds. Thus, by Theorem~B, $F(\vec{z})=\vec{z}$, for
some $\vec{z}\in\Pi$, whence $f_k(\vec{z})=0$, $k=\overline{1,n}$.
Lemma 1 is proved.
\end{proof}
\noindent
{\bf Lemma 2.} Let $f$ be an increasing continuous function on
$[\,0,1\,]$ such that $f(0)=0$, $f(1)=1$ and $f'(0)>\frac 12\cdot
9^n$. Then there exist the numbers $0< z_1<z_2<\ldots<z_n\le 1$
satisfying
\begin{equation}
\label{i9}
f(z_s)=\sum_{k=1}^{n}4\cdot9^{k-1}\frac{z_k^2
z_s}{z_k^{2}+3z_s^{2}},\quad s=\overline{1,n}.
\end{equation}
\begin{proof}
 Since $g(x):=f(x)/x$ is a continuous function on $[\,0,1\,]$ ( $\lim_{x\to 0}g(x)=f'(0)>
\frac 12\cdot 9^n$ ) and $g(1)=1$, then there exist the numbers
$0<u_n<v_n<u_{n-1}<v_{n-1}<\ldots<u_1<v_1\le 1$ for which
$g(u_i)=3\cdot9^{i-1}$ and $g(v_i)=9^{i-1}$, whence $f(u_i)=
3\cdot 9^{i-1}u_i$ and $f(v_i)=9^{i-1}v_i$, $i=\overline{1,n}$.
The fact that $f$ is an increasing function yields
\begin{equation}
\label{i10}
v_i<3\cdot 9^{k-i}u_k,\qquad 1\le k<i\le n.
\end{equation}
For each $s=\overline{1,n}$ put
$$
f_s(\vec{y})=f_s(y_1,\ldots,y_n):=\sum_{k=1}^{n}4\cdot9^{k-1}\frac{y_k^2
y_s}{y_k^{2}+3y_s^{2}}-f(y_s),\quad \vec{y}\in\Pi.
$$
If $\vec{y}\in\Pi_s^+$, then $y_s=v_s$, hence
$$
f_s(\vec{y})> 4\cdot9^{s-1}\frac{v_s^2
v_s}{v_s^{2}+3v_s^{2}}-f(v_s)=0,\quad s=\overline{1,n}.
$$
If $\vec{y}\in\Pi_s^-$, then $y_s=u_s$, hence
\begin{align*}
f_s(\vec{y})&=\sum_{k=1}^n4\cdot9^{k-1}\frac{y_k^2
u_s}{y_k^{2}+3u_s^{2}}-f(u_s)=
\sum_{k=1}^{s-1}4\cdot9^{k-1}\frac{y_k^2
u_s}{y_k^{2}+3u_s^{2}}\\
&+9^{s-1}u_s+\sum_{k=s+1}^{n}4\cdot9^{k-1}\frac{y_k^2
u_s}{y_k^{2}+3u_s^{2}}-3\cdot 9^{s-1}u_s\\
&\le\sum_{k=1}^{s-1}4\cdot9^{k-1}u_s+\sum_{k=s+1}^{n}4\cdot9^{k-1}
\frac{v_k^2}{3u_s^2}u_s-2\cdot 9^{s-1}u_s\\
&\le\frac 12\cdot 9^{s-1}u_s+\sum_{k=s+1}^{n}\frac 43\cdot
9^{2s-k}u_s-2\cdot 9^{s-1}u_s<0,\quad s=\overline{1,n},
\end{align*}
where in the last line we use~\eqref{i10}. Applying Lemma 1 for
the functions $f_s$, $s=\overline{1,n}$ we get that there exists
$\vec{z}=(z_1,\ldots z_n)\in\Pi$, such that $f_s(\vec{z})=0$,
$s=\overline{1,n}$, which is~\eqref{i9}. Lemma 2 is proved.
\end{proof}

\section{Proof of Theorem 1}
Let $R\in Q_{2n}$ be an odd and increasing on function on
$[\,-1,1\,]$ such that $R'(0)>\frac 12\cdot 9^nR(1)$. Without any
loss of generality we may assume that $R(1)=1$. By Lemma 2, there
exist the numbers $0<z_n<z_{n-1}<\ldots<z_1\le 1$ such that the
function
$$
L(x):=\sum_{k=1}^{n}4\cdot9^{k-1}\frac{z_k^2
x}{z_k^{2}+3x^{2}}-R(x)
$$
satisfies the equalities $L(z_s)=0$, $s=\overline{1,n}$. Further,
we have
$$
L'(0)=\sum_{k=1}^{n}4\cdot9^{k-1}-R'(0)<0,
$$
and
\begin{align*}
L'(z_s)&=\sum_{k=1}^{n}4\cdot9^{k-1}\frac{z_k^2(z_k^2-3z_s^2)}{(z_k^2+3z_s^2)^2}
-R'(z_s)\\
&<\sum_{k=1}^{s-1}4\cdot9^{k-1}-\frac 12\cdot 9^{s-1}<0.
\end{align*}
Thus, each of the open intervals $(0, z_n)$, $(z_n, z_{n-1})$,
...,$(z_2, z_1)$ contains at least one zero of the function $L$.
Since $L$ is an odd function, then $L$ has at least $4n+1$ zeroes
on $[\,-1,1\,]$. On the other hand $l\in Q_{4n}$, so $L\equiv 0$.
This contradiction finished the proof of Theorem 1.\\
{\bf Proof of Corollary 1:} Without any loss of generality we may
assume that $x>0$. For each increasing rational function $R\in
Q_n$ and $x\in (0,1)$ put
$$
H(y):=\frac{R(x+y(1-x))-R(x-y(1-x))}2.
$$
Evidently, $H\in Q_{2n}$ is odd increasing rational function with
$\|H\|\le \|R\|$. Note that $H'(0)=(1-x)R'(x)$. Thus, applying
Theorem 1 for the function $H$ we get
$$
R'(x)\le\frac{9^n}{2(1-x)}\|R\|\le\frac{ 9^n}{1-x^2}\|R\|.
$$
\section{Proof of Theorem 2}
Below we use without proof two easy inequalities: If $\gamma$,
$\alpha>0$, then
\begin{equation}
\label{T21}
\frac{\gamma^2(\gamma^2-x^2)}{(\gamma^2+x^2)^2}\ge
-\frac{\gamma^2}{\alpha^2},\quad |x|\ge\alpha,
\end{equation}
and
\begin{equation}
\label{T22}
\frac{\gamma^2(\gamma^2-x^2)}{(\gamma^2+x^2)^2}\ge
-\frac 18,\quad x\in\mathbb{R}.
\end{equation}
Let $K_n$ be the set of all odd rational functions $R\in Q_{2n-1}$
with $R'(x)>0$, $x\in [-1,1]$. Evidently, it is sufficient to
prove that
\begin{equation}
\label{T23}
S_n:=\sup_{R\in K_n}\frac{R'(0)}{\|R\|}\ge 9^{n-1}.
\end{equation}
If $n=1$, then the function $R(x)\equiv x$ provide~\eqref{T23}.
Let $R\in K_n$ be an arbitrary function. Fix $\epsilon\in (0,
R'(0)/2)$. Since $R'$ is a continuous function, then there exists
$\alpha$, $\beta>0$ such that $R'(x)>R'(0)-\epsilon$, for
$|x|<\alpha$, and $R'(x)>\beta$, for $x\in [-1,1]$. For each
$\gamma>0$ put
$$
G_{\gamma}(x):=R(x)+8(R'(0)-2\epsilon)\frac{\gamma^2x}{\gamma^2+x^2},
\quad x\in[-1,1].
$$
We have
\begin{equation}
\label{T24}
G'_{\gamma}(x)=R'(x)+8(R'(0)-2\epsilon)\frac{\gamma^2(\gamma^2-x^2)}{(\gamma^2+x^2)^2},
\end{equation}
so
\begin{equation}
\label{T25}
G'_{\gamma}(0)=9R'(0)-16\epsilon.
\end{equation}
The inequality~\eqref{T22} implies
$G'_{\gamma}(x)>R'(0)-\epsilon-(R'(0)-2\epsilon)=\epsilon$, for
$|x|<\alpha$. Thus, by~\eqref{T21} and~\eqref{T24} $G_{\gamma}\in
K_{n+1}$, for all $\gamma$ small enough, such that
$$
8(R'(0)-2\epsilon)\frac{\gamma^2}{\alpha^2}<\beta.
$$
Moreover
\begin{equation}
\label{T26}
\|G_{\gamma}\|\le
\|R\|+4(R'(0)-2\epsilon)\gamma \to \|R\|,
\quad
\gamma\to 0.
\end{equation}
Since $R\in K_n$ is an arbitrary function and $\epsilon$ can be
arbitrary small, then~\eqref{T25} and~\eqref{T26} yield
$S_{n+1}\ge 9S_n$. This gives us~\eqref{T23}. Theorem 2 is proved.


\begin{thebibliography}{xx}
\bibitem{LMG}
G.\ G.\ Lorentz,  M.\ v.\ Golitschek, Y.\ Makovoz,  Constructive
Approximation, Springer Verlag, Berlin, 1996.
\bibitem{N}
 L.\ Nirenberg, Topics in nonlinear functional analysis, New York,
 1974.
\end{thebibliography}
\end{document}